# Semi-Selfdecomposable Laws in the Minimum Scheme

S. Satheesh

*NEELOLPALAM, S. N. Park Road, Trichur – 680 004, **India**.*

*ssatheesh@sancharnet.in*

E. Sandhya

*Department of Statistics, Prajyoti Niketan College, Pudukkad*
*Trichur – 680 301, **India**.*

*esandhya@hotmail.com*

**Abstract**

We discuss semi-selfdecomposable laws in the minimum scheme and characterize them using an autoregressive model. Semi-Pareto and semi-Weibull laws of Pillai (1991) are shown to be semi-selfdecomposable in this scheme. Methods for deriving this class of laws are then attempted from the angle of randomization. Finally, discrete analogues of these results are also considered.

## 1 Introduction

Motivated by the formulation of semi-selfdecomposable (SSD) laws in the classical (additive) scheme and max-SSD laws in the maximum scheme, here we discuss min-SSD laws in the minimum scheme. Certain aspects of SSD laws and related processes have been discussed in the classical and maximum schemes in Satheesh and Sandhya (2004, 2006*b*). Min-SSD laws were introduced in Satheesh, *et al*. (2005) to deal with a generalization of the marginally stationary autoregressive model with a minimum structure. Here we consider these laws in some more detail.

Let $X_1$, $X_2$, …. are i.i.d with d.f $F(x)$ and let $U = \text{Minimum}\{X_1, X_2, …., X_n\}$. Then the d.f of $U$ is $1-\{1-F(x)\}^n$. Thus it is more convenient to work with the survival function (s.f) $S(x)$ of $F(x)$ while discussing the distribution of minimums. Also distributions of the minimum are important in the context of stochastic models modeling a series system in reliability and minification processes in time series.

This investigation of min-SSD laws is motivated by the corresponding results in the classical and maximum structures in Satheesh and Sandhya (2004, 2006*b*) and the possibility of characterizing a minification process in the context of an autoregressive time series. Here we will use the symbol $\wedge$ for minimum.

Let $X_1$, $X_2$, …. be a sequence of r.vs. Then this sequence describes a first order autoregressive model with a minimum structure (min-AR(1)) if for some $\rho>0$ there exists another sequence of i.i.d r.vs $\{\varepsilon_i\}$ such that

$$X_n \stackrel{d}{=} \rho X_{n-1} \wedge \varepsilon_n, \text{ for all } n>0 \text{ integer.} \qquad (1.1)$$

In section.2 we characterize min-SSD laws as the one that models a marginally stationary min-AR(1) scheme and show that the semi-Pareto and semi-Weibull laws are min-SSD. Methods for deriving min-SSD laws are discussed in section.3 using two ways of randomization. In section.4 we describe the integer-valued analogue of these distributions and an integer-valued min-AR(1) model.

## 2  Min-Semi-Selfdecomposable Laws

First we make a general observation regarding a s.f $S(x)$ on **R**.

**Lemma.2.1** If $S(x)$ is a s.f of a continuous r.v $X$ then $\{S(x)\}^a$ is a s.f for any $a>0$.

*Proof*. The characteristic properties of a s.f *viz*. $S(-\infty) = 1$, $S(\infty) = 0$ and $S(x)$ is continuous are intact with $\{S(x)\}^a$ also. Finally, that $\{S(x)\}^a$ is non-increasing is clear by taking its first derivative.

**Definition.2.1** A non-degenerate d.f $F$ with s.f $S$ is min-SSD($b$) if there exists another s.f $S_o$ such that

$$S(x) = S(bx) S_o(x), \text{ for all } x \in \mathbf{R} \text{ and for some } b \in (0,1) \cup (1,\infty). \qquad (2.1)$$

If this is true for every $b \in (0,1) \cup (1,\infty)$, then $F$ is min-SD.

Let us now consider the min-AR(1) model (1.1). Since $X_{n-1}$ is a function only of $\varepsilon_j$, $j=1, 2, ...., n-1$, it is independent of $\varepsilon_n$. Hence in terms of s.fs this is equivalent to;

$$S_n(x) = S_{n-1}(x/\rho) S_\varepsilon(x), \text{ for all } x \in \mathbf{R}.$$

Assuming the series to be marginally stationary we have;

$$S(x) = S(x/\rho)\, S_\varepsilon(x). \tag{2.2}$$

Now comparing (2.2) and (2.1) the following theorem is clear.

**Theorem.2.1** A sequence $\{X_n\}$ of r.vs generates the min-AR(1) series (1.1) that is marginally stationary iff the distribution of $X_n$ is min-SSD($1/\rho$).

Semi-Pareto and semi-Weibull are two models that are discussed in the minification structure in Pillai (1991). We now show that these laws are min-SSD.

**Definition.2.2** Semi-Pareto($p,\alpha$) family of laws are those with s.f $1/(1+\psi(x))$, $x>0$, where $\psi(x)$ satisfies $p\psi(x) = \psi(p^{1/\alpha}x)$, for all $x>0$, some $0<p<1$, and $\alpha>0$.

**Definition.2.3** Semi-Weibull($p,\alpha$) family of laws are those with s.f $\exp\{-\psi(x)\}$, $x>0$, where $\psi(x)$ satisfies $p\psi(x) = \psi(p^{1/\alpha}x)$, for all $x>0$, some $0<p<1$, and $\alpha>0$.

**Theorem.2.2** Semi-Pareto($p,\alpha$) family of laws is min-SSD($p^{1/\alpha}$).

*Proof.* Pillai and Sandhya (1996) showed that the geometric minimum of a semi-Pareto($p,\alpha$) law is of its own type (the property of geometric-min semi-stability), the geometric law being on $I_1=\{1, 2, \ldots\}$. That is, its s.f satisfies

$$\frac{1}{1+\psi(x)} = \frac{p\frac{1}{1+\psi(cx)}}{1-(1-p)\frac{1}{1+\psi(cx)}} = \frac{1}{1+\frac{1}{p}\psi(cx)}, \text{ for all } x>0 \text{ and some } c>0.$$

Hence, $\psi(x) = \frac{1}{p}\psi(cx)$, where the constants $\frac{1}{p}$ and $c$ are related by $\frac{1}{p}c^\alpha = 1$, and hence $c = p^{1/\alpha}$. Thus $\psi(x)$ satisfies the condition in definition.2.2. Now, the above equation can also be written as:

$$\frac{1}{1+\psi(x)} = \frac{1}{1+\psi(cx)}\frac{p}{1-(1-p)\frac{1}{1+\psi(cx)}}.$$

Here the second factor on the RHS is also a s.f being the geometric minimum of the semi-Parto law where the geometric law is on $I_o=\{0, 1, 2, \ldots\}$. Hence semi-Pareto family is min-SSD($p^{1/\alpha}$).

**Theorem.2.3** Semi-Weibull$(p,\alpha)$ family of laws is min-SSD$(p^{1/\alpha})$.

*Proof.* We have the s.f $exp\{-\psi(x)\} = exp\{-a\psi(p^{1/\alpha}x)\}$, $a = 1/p$.

$$= exp\{-\psi(p^{1/\alpha}x)\} \, exp\{-(a-1)\psi(p^{1/\alpha}x)\}.$$

Here both the factors on the RHS are s.fs, the second factor by lemma.2.1, and hence semi-Weibull$(p,\alpha)$ family is min-SSD$(p^{1/\alpha})$.

In the next section we extend the main line of argument in theorem.2.2 to discuss a method to derive min-SSD laws and then extend it in another direction.

### 3  Methods for Deriving Min-SSD Laws

Notice that in theorem.2.2 the main idea was to write the s.f as two factors as a consequence of the geometric minimum semi-stability of semi-Pareto laws. Again this was possible since the geometric law was supported by $I_1$. To generalize this we need the following notion of N-min (semi-)stability, see Satheesh and Nair (2002).

**Definition.3.1** Let $X_1, X_2, \ldots$ be non-degenerate i.i.d r.vs with a common d.f $F$ (s.f $S$) and $N$ be a positive integer-valued r.v that is independent of $X$ with probability generating function (PGF) $Q(s)$. Then $F$ is N-min semi-stable if:

$$Q\{S(cx)\} = S(x) \text{ for all } x\in \mathbf{R} \text{ and some } c>0. \tag{3.1}$$

If this relation is true for all $c>0$ then $F$ is N-min stable.

**Remark.3.1** When $N$ is degenerate at $k>0$ integer then $F$ is semi-Weibull$(1/k,\alpha)$, see Satheesh and Nair (2002).

**Theorem.3.1** Every distribution $F$ that is N-min semi-stable is min-SSD.

*Proof.* Let $Q(s)$ be the PGF of the r.v $N$ that is positive and integer-valued and $n$ be the starting point of the support of $N$. Hence $N \stackrel{d}{=} 1+M$, where the starting point of the support of the integer-valued r.v $M$ is $(n\text{-}1)$. That is, the starting point of the support of $M$ is at least zero and let us denote its PGF by $P(s)$. Hence $Q(s) = sP(s)$. Now, since $F$ is N-min semi-stable $Q\{S(cx)\} = S(x)$ for all $x\in \mathbf{R}$ and some $c>0$. That is;

$$S(x) = Q\{S(cx)\} = S(cx) \, P\{S(cx)\},$$

where $P\{S(cx)\}$ is another s.f. Hence $F$ is min-SSD$(c)$.

**Corollary.3.1** Every distribution that is N-min stable is min-SD.

Satheesh and Sandhya (2006*a, b*) have introduced φ-max-semi-stable laws as: For a Laplace transform (LT) $\varphi$ the d.f $\varphi\{-ln(F(x))\}$ is φ-max-semi-stable if $F(x)$ is max-semi-stable. Similarly here we consider φ-semi-Weibull laws as those with s.f $\varphi\{\psi(x)\}$ where $\varphi$ is a LT and $exp\{-\psi(x)\}$ is the s.f of the semi-Weibull($p,\alpha$) laws. Here essentially we are randomizing the parameter *a* in $\{S(x)\}^a$ (lemma.2.1) for the semi-Weibull family. Notice that semi-Pareto laws are exponential-semi-Weibull laws and we know that semi-Weibull laws are min-SSD and the exponential law is SD. We now describe φ-semi-Weibull laws that are min-SSD. We need the following notion.

**Definition.3.2** (Maejima an Naito, 1998). A LT $\varphi$ is SSD(*c*) if for some $0<c<1$ there exists a LT $\varphi_o(s)$ such that

$$\varphi(s) = \varphi(cs)\,\varphi_o(s), \text{ for all } s>0.$$

If this relation holds for every $c \in (0,1)$ then the distribution with LT $\varphi$ is SD.

**Theorem.3.2** φ-semi-Weibull($p,\alpha$) laws are min-SSD($p^{1/\alpha}$) if $\varphi$ is SSD($p$).

*Proof.* We have the s.f $\varphi\{\psi(x)\}$ where $\varphi$ is a LT and $p\psi(x) = \psi(p^{1/\alpha}x)$, for all $x>0$, some $0<p<1$, and $\alpha>0$. If $\varphi$ is SSD(*p*) then there exists another LT $\varphi_o(s)$ and,

$$\varphi(s) = \varphi(ps)\,\varphi_o(s), \text{ for all } s>0 \text{ and some } 0<p<1.$$

Hence the s.f of the corresponding φ-semi-Weibull law can be written as;

$$\varphi\{\psi(x)\} = \varphi\{p\psi(x)\}\,\varphi_o\{\psi(x)\} = \varphi\{p\frac{1}{p}\psi(p^{1/\alpha}x)\}\,\varphi_o\{\psi(x)\}$$

$$= \varphi\{\psi(p^{1/\alpha}x)\}\,\varphi_o\{\psi(x)\}, \text{ completing the proof.}$$

**Corollary.3.2** φ-semi-Weibull($p,\alpha$) laws are min-SSD($p^{1/\alpha}$) if $\varphi$ is SD.

**Example.3.1** A restatement of proposition.3 in Satheesh and Nair (2002) is: The non-negative i.i.d r.vs $X_i$ are Harris(1,*a,k*)-min semi-stable iff they are generalized semi-Pareto($p,\alpha,1/k$) with d.f $F(x) = 1-\{1+\psi(x)\}^{-1/k}$, $x>0$, $p\psi(x) = \psi(p^{1/\alpha}x)$, for all $x>0$, some $0<p<1$, $p=1/a$, $k \in I_1$ and $\alpha>0$. The PGF of the Harris(1,*a,k*) law is $s/\{a-(a-1)s^k\}^{1/k}$. This being the PGF of a positive integer-valued r.v $\{a-(a-1)s^k\}^{-1/k}$ is also a PGF. Hence by theorem.3.1 the generalized semi-Pareto($p,\alpha,1/k$) family of laws are min-SSD.

**Example.3.2** Pareto laws with s.f $1/\{1+x^{\alpha}\}$, $x>0$, $\alpha>0$, are min-SD as they are geometric-min stable, the geometric law being on $I_1$.

**Example.3.3** Consider the s.f $S(x) = \{1+\psi(x)\}^{-\beta}$, $x>0$, $\beta>0$ and $p\psi(x) = \psi(p^{1/\alpha}x)$, for all $x>0$, some $0<p<1$ and $\alpha>0$, of the generalized semi-Pareto$(p,\alpha,\beta)$ family of laws. Since the gamma$(1,\beta)$ law is SD the above family is min-SSD$(p^{1/\alpha})$ being that of a gamma-semi-Weibull law.

**Remark.3.2** it may be noted that the above conclusion cannot be arrived at using the approach in example.3.1 because the Harris-min-semi-stability implies that the distribution is generalized semi-Pareto$(p,\alpha,1/k)$ only and not generalized semi-Pareto$(p,\alpha,\beta)$, see Satheesh and Nair (2002).

**Example.3.4** By a line of argument similar to theorem.3.2, Satheesh and Sandhya (2004) have shown that the characteristic function $\{1+\psi(t)\}^{-\beta}$ is SSD($b$), $b<1$ (Hence the LT $\{1+\psi(s)\}^{-\beta}$ is also SSD($b$)). Here $\psi$ is such that $\psi(s) = a\psi(bs)$, $\forall s>0$ and some $0<b<1<a$, $ab^{\alpha}=1$, $0<\alpha<1$. Setting $\varphi$ to be the LT of this SSD($b$) law the corresponding $\varphi$-semi-Weibull laws are min-SSD($b^{1/\alpha}$).

## 4 Integer-Valued Min-AR(1) Schemes

The following construction of integer-valued distributions is from Satheesh and Sandhya (1997). Let $\{m(j)\}$ is the realizations of a LT at $j \in I_0 = \{0,1,2, ....\}$. Now,

$$P\{X<j\} = F(j) = 1- m(j), \quad j \in I_0 \qquad (4.1)$$

is the d.f of a mixture of geometric laws on $I_0$. As $m(cs)$, $c>0$, is again a LT the function $G(j) = 1- m(cj)$, is a d.f. The restriction $\alpha<1$, in the definitions of semi-Weibull, semi-Pareto and generalized semi-Pareto laws makes them mixtures of exponentials (since now $S(x)$'s are LTs) and $x$ to $j \in I_0$, yield their discrete versions as in (4.1), see Satheesh and Nair (2002). Since $F(cj) = G(j)$, definitions of min-SSD and N-min semi-stable laws for d.fs on $I_0$ with expressions analogous to (2.1) and (3.1) holding for all $j \in I_0$ are possible. Similarly the notion of an integer-valued min-AR(1) series on $I_0$ can be described as in (1.1). Thus we have results that are integer-valued analogues of the results in sections 2 and 3. We record them here without proof. These descriptions hold good for distributions on $I_0$ and under the formulation in (4.1).

**Theorem.4.1** A sequence $\{X_n\}$ of non-negative integer-valued r.vs generates a marginally stationary min-AR(1) series if the distribution of $X_n$ is discrete min-SSD($1/\rho$).

**Theorem.4.2** Discrete semi-Pareto($p,\alpha$) and Discrete semi-Weibull($p,\alpha$) family of laws are min-SSD($p^{1/\alpha}$).

The s.f of a semi-Weibull($p,\alpha$) law with $\alpha<1$, is the LT of a semi-stable law which is infinitely divisible (Pillai, 1971). Hence any real power of this LT is again a LT and so we can proceed as in the proof of theorem.2.3.

**Theorem.4.3** Every distribution on $I_0$ that is N-min semi-stable is min-SSD.

**Theorem.4.4** Discrete φ-semi-Weibull($p,\alpha$) laws are min-SSD($p^{1/\alpha}$) if $\varphi$ is SSD($p$).

Thus we have some probability laws to model the min-AR(1) series.

*Acknowledgements.* It is our pleasure to place on record the inspiration we received in our research from Professor T S K Moothathu who is to retire in June 2006.

### References.